\documentclass[a4paper,11pt]{article}
\usepackage{amsmath,amscd,amsfonts,amssymb,epsf,latexsym}

\makeatletter

\long\def\@makefnt#1{\parindent 1em\noindent
            \hb@xt@1.8em{\hss\@textsuperscript{}}#1}
\long\def\@ftntext#1{\insert\footins{%
    \reset@font\footnotesize
    \interlinepenalty\interfootnotelinepenalty
    \splittopskip\footnotesep
    \splitmaxdepth \dp\strutbox \floatingpenalty \@MM
    \hsize\columnwidth \@parboxrestore
    \color@begingroup
      \@makefnt{%
        \rule\z@\footnotesep\ignorespaces#1\@finalstrut\strutbox}%
    \color@endgroup}}%
\def\subjclass#1{%
  \@ftntext{2000 {\itshape Mathematics Subject Classification.}\enspace #1.}}
\def\keywords#1{%
  \@ftntext{{\itshape Key words and phrases.}\enspace #1.}}
\makeatother

\def\A{{\mathbb A}}
\def\B{{\mathbb B}}

\def\C{{\mathbb C}}
\def\D{{\mathbb D}}

\def\X{{\mathbb X}}
\def\Y{{\mathbb Y}}
\def\E{{\mathbb E}}

\def\AB{ {\mathbb A}\moins {\mathbb B}}

\def\moins{\raise 1pt\hbox{{$\scriptstyle -$}}}
\def\plus{\raise 1pt\hbox{{$\scriptstyle +$}} }

\def\phi{\varphi}

\newtheorem{theorem}{Theorem}
\newtheorem{proposition}[theorem]{Proposition}
\newtheorem{lemma}[theorem]{Lemma}
\newtheorem{corollary}[theorem]{Corollary}
\newtheorem{remark}[theorem]{Remark}
\newtheorem{definition}[theorem]{Definition}

\newtheorem{convention}[theorem]{Convention}
\newtheorem{example}[theorem]{Example}

\def\proof{\noindent{\bf Proof.\ }}

\def\qed{~\hbox{$\Box$}}

\def\Aut{\mathop{\rm Aut}}
\def\Diff{\mathop{\rm Diff}}
\def\codim{\mathop{\rm codim}}

\def\rank{\mathop{\rm rank}}

\begin{document}

\title{\bf Thom polynomials and Schur functions I}

\author{Piotr Pragacz\thanks{Research partially supported by KBN
grant 2P03A 024 23. Part of the work was done during the author's stay
at the METU in Ankara (May 2005), supported by T\"UB\.ITAK.}\\
\small Institute of Mathematics of Polish Academy of Sciences\\
\small \'Sniadeckich 8, 00-956 Warszawa, Poland\\
\small P.Pragacz@impan.gov.pl}

\date{(30.08.2005; revised 15.05.2006)}



\maketitle

\centerline{\it To the memory of Professor Stanis\l aw Balcerzyk
(1932-2005)}

\begin{abstract}
We give the Thom polynomials -- via their Schur function
expansions -- for the singularities $I_{2,2}$,
and $A_3$ associated with maps $({\bf C}^{\bullet},0)
\to ({\bf C}^{\bullet+k},0)$ with parameter $k\ge 0$.
Moreover, for the singularities $A_i$ (with any parameter $k\ge 0$)
we analyze the ``first approximation'' $F^{(i)}$ to the Thom polynomial.
Our computations combine the characterization of Thom polynomials
via the ``method of restriction equations'' of Rimanyi et al. with
the techniques of (super) Schur functions.
\end{abstract}

\tableofcontents

\section{Introduction}\label{intro}

The global behavior of singularities is governed by their {\it Thom
polynomials} (cf. \cite{T}, \cite{Kl}, \cite{AVGL}, \cite{Ka},
\cite{Rim2}). Knowing the Thom polynomial of a singularity $\eta$,
denoted ${\cal T}^{\eta}$, one can compute the cohomology class represented
by the $\eta$-points of a map.
We do not attempt here to survey all activities related
to computations of Thom polynomials -- note that in
Kleiman's survey \cite{Kl},
the reader can find a summary of an early work of Thom, Porteous, Ronga,
Menn, Sergeraert, Lascoux, and Roberts on Thom polynomials presented
in an algebro-geometric framework.

In the present paper, following a series of papers by Rimanyi et al.
\cite{RS}, \cite{Rim1}, \cite{Rim2}, \cite{FR}, \cite{BFR},
we study the Thom polynomials for the singularities
$I_{2,2}$ and $A_i$ associated with maps
$({\bf C}^{\bullet},0) \to ({\bf C}^{\bullet+k},0)$
with parameter $k\ge 0$.

The way of obtaining the thought Thom polynomial is through the solution
of a system of linear equations, which is fine when we want to find
one concrete Thom polynomial, say, for a fixed $k$. However, if we want
to find the Thom polynomials for a series of singularities, associated
with maps $({\bf C}^{\bullet},0) \to ({\bf C}^{\bullet+k},0)$
with $k$ as a parameter, we have to solve {\it simultaneously} a countable
family of systems of linear equations. As stated by Rimanyi in \cite{Rim2},
p. 512~:

\smallskip

{``However, another challenge is to find Thom polynomials containing
$k$ as a parameter.''}

\smallskip

We do it here for the restriction
equations for the above mentioned singularities. In fact, the obtained
functional equations in symmetric functions are of independent interest.
The main novelty of
the present paper over the previous articles on Thom polynomials,
is an extensive use of {\it Schur functions}. Namely, instead of using
{\it Chern monomial expansions} (as the authors of all previous papers
constantly do), we use {\it Schur function expansions}.
This puts a more transparent structure on computations of Thom
polynomials. We hope that an expression for ${\cal T}^{I_{2,2}}$
(Theorem \ref{TI22}), as well as our expression for ${\cal T}^{A_3}$
(Theorem \ref{TA3}), provide a support of this claim. For example,
we get in this way some recursive formulas (cf., e.g., Lemma \ref{Lr})
that are not so easy to find using other bases, in particular
the Chern monomial basis. In fact, recursions play a prominent role
in the formulas of the present paper, cf. Eqs. (\ref{rd}) and (\ref{re}).

Another feature of using the Schur function expansions for Thom
polynomials is that in all known to us cases (not only those treated
in the present paper), all the coefficients are nonnegative.\footnote{Added
15.05.2006: Positivity of Schur function expansions of Thom polynomials
has been recently proved by A. Weber and the author in \cite{PW}.}

\smallskip

To be more precise, we use here (the specializations of)
{\it supersymmetric} Schur functions, also called {\it super-S-functions}
or Schur functions {\it in difference of alphabets} together with their
three basic properties: {\it vanishing}, {\it cancellation} and {\it
factorization}, (cf. \cite{BR}, \cite{LS}, \cite{P2}, \cite {PT},
\cite{M}, \cite{FP}, and \cite{L}).
These functions contain {\it resultants} among themselves. Their
geometric significance was illuminated in the 80's in the
author's study of {\it polynomials supported on degeneracy loci}
(cf. \cite{P}). In fact, in the present paper we use the point of view
of that article to some extent.
More precisely, given a morphism $F \to E$ of vector bundles,
where $\rank(E)=m$ and $\rank(F)=n$,
by a $j$-{\it polynomial} we understand a {\bf Z}-linear combination of
the Schur functions \ $S_I(E-F)$, where partitions $I$ are such that
$(n\moins j)^{m-j} \subset I$
but $(n\moins j\plus 1)^{m-j+1}\not\subset I$
\ \footnote{We shall use the same name
for (formal) alphabets $\A$ and $\B$ of cardinalities $m$ and $n$
instead of the alphabets of the
Chern roots of $E$ and $F$.}. In some sense,
a $j$-polynomial is a ``typical'' polynomial supported on the
$j$th degeneracy locus $D_j$ of the morphism
(in the sense of \cite{P}, see also \cite{FP}).

A Thom polynomial is a sum of of such $j$-polynomials associated with the
corresponding morphism of tangent bundles (cf. the next section).
The polynomial ${\cal T}^{I_{2,2}}$ is a {\it single} $j$-polynomial
whereas ${\cal T}^{A_3}$ for $k>0$ is the sum of two $j$-polynomials
(for two consecutive $j$'s).
We first determine its part related to the smaller rectangle,
and then add necessary ``corrections'' related to the larger rectangle.
For the singularities $A_i$ (any $i$), we describe the $j$-polynomial part
of the Thom polynomial for the largest possible $j$, a sort of
the ``first approximation'' to ${\cal T}^{A_i}$ (here, the corresponding
rectangle is a single row of length $k+1$).

\smallskip

We give here a complete proof of our formula for the Thom polynomial
for the singularity $I_{2,2}$. As for the singularities $A_i$,
we reprove the formulas of Thom and Ronga for $A_1$ and $A_2$, and
announce the Schur function expression for $A_3$. We outline a proof
in the Appendix. Another expression for $A_3$ in terms of monomials
in the Chern classes was announced by Berczi, Feher, and Rimanyi
in \cite{BFR}.
We also illustrate the methods used by reproving the
result of Gaffney giving the Thom polynomial for $A_4$ with $k=0$
(this was also done by Rimanyi \cite{Rim1} -- our approach uses more
extensively Schur functions). For any singularity $A_i$, we describe
the $j$-polynomial part (denoted by \ $F^{(i)}_r$) of the Thom
polynomial for the largest possible $j$.

In our calculations we use extensively the $\lambda$-ring approach
to symmetric functions developed mainly in Lascoux's book \cite{L}.

\section{Recollections on Thom polynomials}

Our main reference for this section is \cite{Rim2}.
We start with recalling what we shall mean by a ``singularity''.
Let $k\ge 0$ be a fixed integer. By {\it singularity}
we shall mean an equivalence class of stable germs $({\bf C}^{\bullet},0)
\to ({\bf C}^{\bullet+k},0)$, where $\bullet\in {\bf N}$, under the
equivalence
generated by right-left equivalence (i.e. analytic reparametrizations
of the source and target) and suspension (by suspension of a germ $\kappa$
we mean its trivial unfolding: $(x,v) \mapsto (\kappa(x),v)$).

We recall that the {\it Thom polynomial} ${\cal T}^{\eta}$ of a singularity
$\eta$ is a polynomial in the formal variables $c_1, c_2, \ldots$ that
after the substitution
\begin{equation}
c_i=c_i(f^*TY-TX)=[c(f^*TY)/c(TX)]_i\,,
\end{equation}
for a general map $f:X \to Y$ between complex analytic manifolds,
evaluates the Poincar\'e dual $[\overline{\eta(f)}]$ of the cycle carried
by the closure of
\begin{equation}
\eta(f)=\{x\in X : \hbox{the singularity of} \ f \ \hbox{at} \ x \
\hbox{is} \ \eta \}\,.
\end{equation}
By {\it codimension of a singularity} $\eta$, $\codim(\eta)$,
we shall mean $\codim_X(\eta(f))$ for such an $f$. The concept of
the polynomial ${\cal T}^{\eta}$ comes from Thom's fundamental paper
\cite{T}.
For a detailed discussion of the {\it existence} of Thom polynomials,
see, e.g., \cite{AVGL}. Thom polynomials associated with group actions
were studied in \cite{Ka}.
In fact, the above is the ``usual case'' with singularities in the region
where moduli (continuous families) of singularities do not occur. This will
be the case of the singularities studied in the present paper. Indeed,
the codimension of all these singularities does not exceed $6k+8$, the
lowest codimension when moduli of singularities start.

According to Mather's classification, singularities are in one-to-one
correspondence with finite dimensional ${\bf C}$-algebras.
We shall use the following notation:

\medskip

-- $A_i$ \ (of Thom-Boardman type $\Sigma^{1_i}$) will stand for the stable
germs with local algebra ${\bf C}[[x]]/(x^{i+1})$, $i\ge 0$;

\medskip

-- $I_{a,b}$ \ (of Thom-Boardman type $\Sigma^2$) for stable germs with
local algebra ${\bf C}[[x,y]]/(xy, x^a+y^b)$, \ $b\ge a\ge 2$;

\medskip

-- $III_{a,b}$ \ (of Thom-Boardman type $\Sigma^2$) for stable germs
with local algebra ${\bf C}[[x,y]]/(xy, x^a, y^b)$, \ $b\ge a\ge 2$
(here $k\ge 1$).

\medskip

Our computations of Thom polynomials for some of the above
singularities, shall use the method which stems from a sequence
of papers by Rimanyi et al. \cite{RS}, \cite{Rim1}, \cite{Rim2},
\cite {FR}, \cite{BFR}.
We sketch briefly this method, refering the interested reader for more
details to these papers, the main references being the last three
mentioned items.

Let $k\ge 0$ be a fixed integer, and let $\eta: ({\bf C}^{\bullet},0) \to
({\bf C}^{\bullet+k},0)$ be a stable singularity with a prototype
$\kappa: ({\bf C}^n,0) \to ({\bf C}^{n+k},0)$. The {\it maximal compact
subgroup of the right-left symmetry group}
\begin{equation}
\Aut \kappa = \{(\phi,\psi) \in \Diff({\bf C}^n,0) \times
\Diff({\bf C}^{n+k},0) : \psi \circ \kappa \circ \phi^{-1} = \kappa \}
\end{equation}
of $\kappa$ will be denoted by $G_\eta$.
Even if $\Aut \kappa$ is much too large to be a finite dimensional
Lie group, the concept of its maximal compact subgroup (up to conjugacy)
can be defined in a sensible way (cf. \cite{Rim2}).
It is clear that $G_\eta$ can be chosen so that images of its projections
to the factors $\Diff({\bf C}^n,0)$ and $\Diff({\bf C}^{n+k},0)$ are
linear. Its representations via the projections on the source ${\bf C}^n$
and the target ${\bf C}^{n+k}$ will be denoted by $\lambda_1(\eta)$ and
$\lambda_2(\eta)$.
The vector bundles associated with the universal principal
$G_\eta$-bundle $EG_\eta \to BG_\eta$ using the representations
$\lambda_1(\eta)$ and $\lambda_2(\eta)$ will be called
$E_{\eta}'$ and $E_{\eta}$. The {\it total Chern
class of the singularity} $\eta$ is defined
in $H^{\bullet}(BG_\eta;{\bf Z})$ by
\begin{equation}
c(\eta):=\frac{c(E_{\eta})}{c(E_{\eta}')}\,.
\end{equation}
The {\it Euler class} of $\eta$ is defined in
$H^{2\codim(\eta)}(BG_\eta;{\bf Z})$ by
\begin{equation}
e(\eta):=e(E_{\eta}')\,.
\end{equation}

In the following theorem we collect the information from \cite{Rim2},
Theorem 2.4 and \cite{FR}, Theorem 3.5, needed for the calculations
in the present paper.

\begin{theorem}\label{TEq} \ Suppose, for a singularity $\eta$, that
the Euler classes of all singularities of smaller codimension than
$\codim(\eta)$, are not zero-divisors \footnote{This is the so-called
``Euler condition'' (loc.cit.). The Euler condition holds true for the
singularities in the present paper.}. Then we have

\noindent
(i) \ if \ $\xi\ne \eta$ \ and \ $\codim(\xi)\le \codim(\eta)$, then
\ ${\cal T}^{\eta}(c(\xi))=0$;

\noindent
(ii) \ ${\cal T}^{\eta}(c(\eta))=e(\eta)$.

\noindent
This system of equations (taken for all such $\xi$'s) determines the
Thom polynomial ${\cal T}^{\eta}$ in a unique way.
\end{theorem}

To use this method of determining the Thom polynomials for
singularities, one needs their classification, see, e.g., \cite{dPW}.

\medskip

In the forthcoming sections, we shall use these equations to compute
Thom polynomials. Sometimes it will be convenient not to work with
the whole maximal compact subgroup $G_\eta$ but with its suitable subgroup;
this subgroup should be, however, as ``close'' to $G_\eta$ as possible
(cf. \cite{Rim2}, p. 502).

Being challenged by \cite{Rim2}, p. 512 and especially \cite{BFR},
we shall find Thom polynomials {\it containing $k$ as a parameter}
-- this seems to be a (much) more difficult task than computing Thom
polynomials for separate values of $k$, because one must solve
{\it simultaneously} a countable family of systems of linear equations.

\medskip

To effectively use Theorem \ref{TEq} we need to study the maximal
compact subgroups of singularities. We recall the following recipe
from \cite{Rim2} pp. 505--507. Let $\eta$ be a singularity whose
prototype is $\kappa: ({\bf C}^n,0)\to ({\bf C}^{n+k},0)$. The germ
$\kappa$ is the miniversal unfolding of another germ $\beta:
({\bf C}^m,0)\to ({\bf C}^{m+k},0)$ with $d\beta=0$. The group $G_\eta$
is a subgroup of the maximal compact subgroup of
the algebraic automorphism group of
the local algebra $Q_\eta$
of $\eta$ times the unitary group $U(k\moins d)$, where $d$
is the difference between the minimal number of relations and the number
of generators of $Q_\eta$.
With $\beta$ well chosen, $G_\eta$ acts as right-left symmetry group
on $\beta$ with representations $\mu_1$ and $\mu_2$. The representations
$\lambda_1$ and $\lambda_2$ are
\begin{equation}
\lambda_1=\mu_1\oplus \mu_V \ \ \hbox{and} \ \
\lambda_2=\mu_2\oplus \mu_V\,,
\end{equation}
where $\mu_V$ is the representation of $G_\eta$ on the unfolding space
$V={\bf C}^{n-m}$ given, for $\alpha \in V$ and $(\phi,\psi)\in G_\eta$,
by
\begin{equation}
(\phi,\psi) \ \alpha = \psi \circ \alpha \circ \psi^{-1}\,.
\end{equation}
For example, for the singularity of type $A_i$: $({\bf C}^{\bullet},0) \to
({\bf C}^{\bullet+k},0)$, we have $G_{A_i}=U(1)\times U(k)$ with
\begin{equation}
\mu_1=\rho_1, \ \ \mu_2=\rho_1^{i+1}\oplus \rho_k, \ \
\mu_V=\oplus_{j=2}^i \ \rho_1^j \oplus \oplus_{j=1}^i (\rho_k \otimes
\rho_1^{-1})\,,
\end{equation}
where $\rho_j$ denotes the standard representation of
the unitary group $U(j)$. Hence we obtain assertion (i) of the following

\begin{proposition}\label{Pce}
(i) \ Let $\eta=A_i$; for any $k$, writing $x$ and $y_1$,\ldots, $y_k$
for the Chern roots of the universal bundles on $BU(1)$ and $BU(k)$,
\begin{equation}
c(A_i)=\frac{1+(i+1)x}{1+x}\prod_{j=1}^k(1+y_j)\,,
\end{equation}
\begin{equation}\label{eA}
e(A_i)= i! \ x^i \ \prod_{j=1}^k (ix-y_j)\cdots (2x-y_j)(x-y_j)\,.
\end{equation}

\smallskip

\noindent
(ii) \ Let $\eta=I_{2,2}$; for $k\ge 0$, $G_\eta=U(1)\times U(1)\times
U(k)$,
and writing $x_1, x_2$ for the Chern roots of the universal bundles
on two copies of $BU(1)$ and on $BU(k)$,
\begin{equation}
c(I_{2,2})=\frac{(1+2x_1)(1+2x_2)}{(1+x_1)(1+x_2)}\prod_{j=1}^k(1+y_j)\,,
\end{equation}
\begin{equation}
e(I_{2,2})=x_1x_2(x_1-2x_2)(x_2-2x_1)\prod_{j=1}^k (x_1-y_j)(x_2-y_j)
(x_1+x_2-y_j)\,.
\end{equation}

\smallskip

\noindent
(iii) \ Let $\eta=III_{2,2}$; for $k\ge 1$, $G_\eta=U(2)\times U(k\moins 1)$,
and
writing $x_1, x_2$ and $y_1,\ldots,y_{k-1}$ for the Chern roots of the
universal bundles on $BU(2)$ and $BU(k\moins 1)$,
\begin{equation}
c(III_{2,2})=\frac{(1\plus 2x_1)(1\plus 2x_2)(1\plus x_1\plus x_2)}
{(1\plus x_1)(1\plus x_2)}\prod_{j=1}^{k-1}(1+y_j)\,,
\end{equation}
\begin{equation}\label{eIII}
e(III_{2,2})=(x_1x_2)^2(x_1\moins 2x_2)(x_2\moins 2x_1) \
\prod_{j=1}^{k-1}(x_1\moins y_j) \ \prod_{j=1}^{k-1}(x_2\moins y_j)\,.
\end{equation}
\end{proposition}
These assertions are obtained, in a standard way, following
the instructions of \cite{Rim2}, Sect. 4.

\bigskip

\noindent
{\bf Notational conventions} \
Rather than the Chern classes
$$
c_i(f^*TY-TX)=[f^*c(TY)/c(TX)]_i \,,
$$
we shall use {\it Segre classes} $S_i$ of the virtual
bundle $TX^*-f^*(TY^*)$,
i.e. complete symmetric functions $S_i(\A -\B)$ for the alphabets
of the {\it Chern roots} $\A, \B$ of $TX^*$ and $TY^*$. The reader will
find in the next section a summary of algebraic properties of
the functions $S_i(\A-\B)$, or, more generally, Schur functions
$S_I(\A-\B)$ ($I$ runs over sequences of integers) widely used
in the present paper.

Moreover, it will be more handy to use, instead of $k$, a ``shifted''
parameter
\begin{equation}
r:=k+1\,.
\end{equation}
Sometimes, we shall write the Thom polynomial as ${\cal T}^{\eta}_r$
to emphasize its dependence on $r$.
So, e.g., in our notation, the Thom polynomial for the singularity $A_1$
with $\codim(A_1)=r$ for $r\ge 1$ (in general, $\codim(A_i)=ri$),
will be~: ${\cal T}^{A_1}={\cal T}^{A_1}_r=S_r$, instead of $c_{k+1}$
as in the papers in References. In general, a Thom polynomial in terms of
the $c_i$'s (in those papers) will be written here as a linear combination
of Schur functions obtained by changing each $c_i$ to $S_i$ and expanding
in the Schur function basis. Another example is, for $r=1$, the Thom
polynomial for $A_2$: $c_1^2+c_2$ rewritten in the present notation
as ${\cal T}^{A_2}={\cal T}^{A_2}_1=S_{11}+2S_2$.

\section{Recollections on Schur functions}

In this section we collect needed notions related to symmetric
functions and prove a useful Lemma \ref{LFR}. We adopt the point of view
of \cite{L} for what concerns symmetric functions. Namely, given a
commutative ring, we treat symmetric functions as operators acting
on the ring. (Here, these commutative rings are mostly ${\bf Z}$-algebras
generated by the Chern roots of the vector bundles from Proposition
\ref{Pce}.)

\begin{definition}\label{alph}
By an {\it alphabet} $\A$, we understand a (finite)
multi-set of elements in a commutative ring.
\end{definition}

For $k\in {\bf N}$, by ``an alphabet $\A_k$'' we shall mean an alphabet
$\A=(a_1,\ldots,a_k)$ \ (of cardinality $k$); ditto for
$\B_k=(b_1,\ldots,b_k)$, $\Y_k=(y_1,\ldots,y_k)$, and
$\X_2=(x_1,x_2)$.

\begin{definition}\label{cf}
Given two alphabets $\A$, $\B$, the {\it complete functions} $S_i(\AB)$
are defined by the generating series (with $z$ an extra variable):

\begin{equation}
\sum S_i(\AB) z^i =\prod_{b\in \B} (1\moins bz)/\prod_{a\in \A}
(1\moins az)\,.
\end{equation}
\end{definition}
So $S_i(\A-\B)$ interpolates between $S_i(\A)$ -- the complete
homogeneous symmetric function of degree $i$ in $\A$
and $S_i(-\B)$ -- the $i$th elementary function in $\B$ times $(-1)^i$.

The notation $\A -\B$ is compatible with the multiplication
of series:
\begin{equation}
\sum S_i(\A - \B)z^i \cdot \sum S_j(\A' - \B')z^j =
\sum S_i\bigl((\A+\A') - (\B+\B')\bigr)z^i\,,
\label{}
\end{equation}
the sum $\A + \A'$ denoting the union of two alphabets $\A$ and $\A'$.

\begin{convention} We shall often identify an alphabet
$\A=\{a_1,\ldots,a_m\}$ with the sum $a_1+\cdots +a_m$ and perform usual
algebraic operations on such elements. For example, $\A b$ will
denote the alphabet $(a_1b,\ldots,a_mb)$.
We will give priority to the
algebraic notation over the set-theoretic one. In fact, in the following,
we shall use mostly alphabets of variables.
\end{convention}

We have $(\A+\C) - (\B+\C) = \A-\B$, and this corresponds to
simplification of the common factor for the rational series:
\begin{equation}\label{Canc}
\sum S_i((\A + \C) - (\B + \C))z^i = \sum S_i(\A-\B)z^i\,.
\end{equation}

\begin{definition}\label{sf}
Given a sequence $I=(i_1, i_2, \ldots, i_k)\in {\bf Z}^k$, and alphabets
$\A$ and $\B$, the {\it Schur function} $S_I(\A \moins \B)$ is
\begin{equation}\label{schur}
S_I(\A \moins \B):= \Bigl|
     S_{i_p+p-q}(\A \moins \B) \Bigr|_{1\leq p,q \leq k}  \ .
\end{equation}
\end{definition}
We shall mostly use the case when $I$ is
a partition $I=(0\le i_1\le \cdots \le i_k)$. In fact,
by permuting the columns we see that any determinant of the form
(\ref{schur}) is either zero or is, up to sign, such a determinant
indexed by a partition.
These functions are often called {\it supersymmetric Schur functions}
or {\it Schur functions in difference of alphabets}. Their properties
were studied, among others, in \cite{BR}, \cite{LS}, \cite{P2},
\cite {PT}, \cite{M}, \cite{FP}, and \cite{L}; in the present paper,
we shall use the notation and conventions from this last item).

For example,
$$
S_{33344}(\AB)
=\begin{vmatrix}
S_3 \ & \ S_4 \ & S_5 \ & \ S_7 \ &  \ S_8 \\
S_2 & S_3 & S_4 & S_6 & S_7 \\
S_1 & S_2 & S_3 & S_5 & S_6 \\
1 & S_1 & S_2 & S_4 & S_5 \\
0 & 1 & S_1 & S_3 & S_4
\end{vmatrix}\,,
$$
where $S_i$ means $S_i(\AB)$.

By Eq. (\ref{Canc}), we get the following {\it cancellation property}:
\begin{equation}
S_I((\A + \C) - (\B + \C))=S_I(\A-\B)\,.
\end{equation}

In the following, we shall identify partitions with their {\it Young
diagrams}, as is customary.

We record the following property ({\it loc.cit.}), justifying
the notational
remark from the end of Section 2; for a partition $I$,
\begin{equation}
S_I(\AB)= (-1)^{|I|}S_J(\B \moins \A)=S_J(\B^* \moins \A^*)\,,
\end{equation}
where $J$ is the conjugate partition of $I$ (i.e. the consecutive
rows of $J$ are equal to the corresponding columns of $I$), and
$\A^*$ denotes the alphabet $\{-a_1,-a_2,\ldots \}$.

Fix two positive integers $m$ and $n$.
We shall say that a partition $I=(0<i_1\le i_2\le \cdots \le i_k)$
{\it is contained in} the $(m,n)$-hook if either $k\le m$, or $k> m$
and $i_{k-m}\le n$.
Pictorially, this means that the Young diagram of $I$ is contained
in the ``tickened" hook:

\smallskip

$$
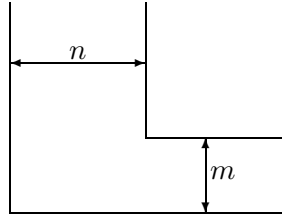
\unitlength=2mm
\begin{picture}(18,14)
\put(0,0){\line(0,1){14}}
\put(0,0){\line(1,0){18}}
\put(9,5){\line(1,0){9}}
\put(9,5){\line(0,1){9}}
\put(4,10){\vector(1,0)5}
\put(4,10){\vector(-1,0)4}
\put(13,2){\vector(0,1)3}
\put(13,2){\vector(0,-1)2}
\put(4.5,10.3){\hbox to0pt{\hss$n$\hss}}
\put(13.3,2.3){\hbox{$m$}}
\end{picture}
$$

\smallskip

We record the following {\it vanishing property}.
Given alphabets $\A$ and $\B$ of cardinalities $m$ and $n$, if
a partition $I$ is not contained in the $(m,n)$-hook,
then ({\it loc.cit.}):
\begin{equation}\label{van}
S_I(\A-\B)=0
\end{equation}
For example,
$$
S_{4569}(\A_2-\B_4)=
S_{4569}(a_1\plus a_2\moins b_1\moins b_2\moins b_3\moins b_4)=0
$$
because $4569$ is not contained in the $(2,4)$-hook.

\smallskip

In fact, we have the following result ({\it loc.cit.}).
\begin{theorem}\label{Tss}
If $\A_m$ and $\B_n$ are alphabets of variables, then the functions
$S_I(\A_m-\B_n)$, for $I$ runing over partitions
contained in the $(m,n)$-hook, are
${\bf Z}$-linearly independent.
\end{theorem}
(They form a ${\bf Z}$-basis of the Abelian group of the so-called
``supersymmetric functions'' ({\it loc.cit.}).)

In the present paper by a {\it symmetric function}, we shall mean
a ${\bf Z}$-linear combination of the operators $S_I(\bullet)$.
By the {\it degree} of such a symmetric function, we shall mean
the largest weight $|I|$ of a partition $I$ involved in its Schur function
expansion.

\medskip
The following useful convention stems from Lascoux's paper \cite{L1}.

\begin{convention} We may need to specialize a letter to $2$, but this must
not be confused with taking two copies of $1$. To allow one, nevertheless,
specializing a letter to an (integer, or even complex) number $r$ inside
a symmetric function, without introducing intermediate variables,
we write \fbox{$r$} for this specialization. Boxes have to be treated
as single variables. For example, $S_i(2) = {{i+1} \choose 2}$ but
$S_i(\fbox{$2$})=2^i$.
A similar remark applies to ${\bf Z}$-linear combinations of variables.
We have \ $S_2(\X_2)=x_1^2\plus x_1x_2\plus x_2^2$ but
$S_2(\fbox{$x_1\plus x_2$})=x_1^2\plus 2x_1x_2\plus x_2^2$,
$S_{11}(\X_2)=x_1x_2$ but
$S_{11}(\fbox{$x_1\plus x_2$})=0$, $S_{2}(3x)=6x^2$ but
$S_{2}(\fbox{$3x$})=9x^2$ etc.
\end{convention}

\begin{definition}
Given two alphabets $\A,\B$, we define their {\it resultant}:
\begin{equation}\label{}
R(\A,\B):=\prod_{a\in \A,\, b\in \B}(a\moins b)\,.
\end{equation}
\end{definition}
This terminology is justified by the fact that $R(\A,\B)$ is the classical
resultant of the polynomials $R(x,\A)$ and $R(x,\B)$.
For example, Eq. (\ref{eA}) can be rewritten as
$$
e(A_i)= R\bigl(x+\fbox{$2x$}+\cdots+\fbox{$ix$},
\Y_k +\fbox{$(i\plus 1)x$} \ \bigr)
$$
and Eq. (\ref{eIII}) is
$$
e(III_{2,2})=R(\X_2,\fbox{$2x_1$}+\fbox{$2x_2$}+\fbox{$x_1\plus x_2$}+
\Y_{k-1})\,.
$$

We have ({\it loc.cit.})
\begin{equation}\label{ER}
R(\A_m,\B_n)= S_{(n^m)}(\AB)=\sum_I S_I(\A) S_{(n^m)/I}(-\B)\,,
\end{equation}
where the sum is over all partitions $I\subset (n^m)$.

When $I$ is contained in the $(m,n)$-hook and at the same time $I$
contains the rectangle $(n^m)$, then we have the following
{\it factorization property} ({\it loc.cit.}):
for partitions $I=(i_1,\ldots,i_m)$ and $J=(j_1,\ldots, j_k)$,
\begin{equation}\label{Fact}
S_{(j_1,\ldots,j_k,i_1+n,\ldots,i_m+n)}(\A_m-\B_n)
=S_I(\A) \ R(\A,\B) \ S_J(-\B)\,.
\end{equation}

\bigskip

We now pass to the following function $F$. Fix positive integers $m$
and $n$. For an alphabet $\A$ of cardinality $m$, we set
\begin{equation}
F(\A, \bullet):= \sum_{I} S_I(\A) S_{n-i_m,\ldots,n-i_1,n+|I|}(\bullet)\,,
\end{equation}
where the sum is over partitions $I=(i_1\le i_2 \le \cdots \le i_m \le n)$.

\begin{lemma}\label{LFR} For a variable $x$ and an alphabet $\B$
of cardinality $n$,
\begin{equation}
F(\A,x-\B)= R(x+\A x,\B)\,.
\end{equation}
\end{lemma}
\proof
For a fixed partition $I=(i_1\le i_2 \le \cdots \le i_m\le n)$, it follows
from the factorization property (\ref{Fact}) that
$$
S_{n-i_m,\ldots,n-i_1,n+|I|}(x-\B)=S_{(n^m)/I}(-\B) \ R(x, \B) \ x^{|I|}\,.
$$
Hence, using $S_I(\A x)= S_I(\A) x^{|I|}$, a standard factorization
of a resultant, and Eq. (\ref{ER}), we have
$$
\aligned
F(\A,x-\B)=\sum_I S_I(\A) S_{(n^m)/I}(-\B) \ R(x, \B) \ x^{|I|}\\
=\sum_I S_I(\A x) \ S_{(n^m)/I}(-\B) \ R(x, \B)& \\
=R(\A x,\B) \ R(x,\B)
=R(x+\A x, \B)& \,.
\endaligned
$$
The lemma has been proved. \qed

We end this section, with an example illustrating the way we shall
use symmetric functions in our computations.
We consider the singularity $III_{2,2}$ (with parameter $r\ge 2$) whose
codimension is $2r+2$. We know that ${\cal T}^{III_{2,2}}$ is equal to the
Thom polynomial for $\Sigma^2$, and the latter polynomial was computed
in \cite{Po}.

But let us apply directly Theorem \ref{TEq}
to the singularity $III_{2,2}$. By virtue of Proposition \ref{Pce},
for $r\ge 3$, the equations characterizing the Thom polynomial for
$III_{2,2}$ are:
\begin{equation}\label{v1}
P(-\B_{r-1})=P(x-\fbox{$2x$}-\B_{r-1})=P(x-\fbox{$3x$}-\B_{r-1})=0\,,
\end{equation}
and additionally,
\begin{equation}\label{norm}
P(\X_2 \moins \fbox{$2x_1$} \moins \fbox{$2x_2$} \moins \fbox{$x_1 \plus x_2$}
\moins \B_{r-2})
=R(\X_2,\fbox{$2x_1$} \plus \fbox{$2x_2$} \plus \fbox{$x_1+x_2$}
\plus \B_{r-2})\,.
\end{equation}
Here, without loss of generality, we assume that $x$, $x_1, x_2$, and
$\B_{r-1}$ are variables, and $P(\bullet)$
denotes a symmetric function. Indeed, the singularities $\ne III_{2,2}$
of codimension $\le \codim(III_{2,2})$ are: $A_0$, $A_1$, $A_2$.
For $r=2$ we must add the vanishing imposed by $A_3$ which
(similarly to $III_{2,2}$) is of codimension $6$ (this is the only
exception):
\begin{equation}\label{v2}
P(x-\fbox{$4x$}-\B_1)=0\,.
\end{equation}
Since the partition $(r \plus 1,r \plus 1)$ is not contained in the
$(1,r)$-hook, for $P=S_{r+1,r+1}(\bullet)$ we get the vanishings
(\ref{v1}) and (\ref{v2}).
Moreover, Eq. (\ref{norm}) is satisfied for this $P$ because
$$
P(\X_2 \moins \fbox{$2x_1$} \moins \fbox{$2x_2$} \moins
\fbox{$x_1+x_2$} \moins \B_{r-2})=
R(\X_2,\fbox{$2x_1$} \plus \fbox{$2x_2$} \plus \fbox{$x_1 \plus x_2$}
\plus\B_{r-2})\,.
$$
These equations characterize the Thom polynomial for $III_{2,2}$,
and hence this polynomial is equal to $S_{r+1,r+1}$ in agreement with
\cite{Po}. In the forthcoming computations, however, the method
of restriction equations from Theorem \ref{TEq} will play a principal role.

\section{Thom polynomial for $I_{2,2}$}

The codimension of $I_{2,2}$ (for parameter $r\ge 1$) is $3r+1$.
For $r=1$, the Thom polynomial for $I_{2,2}$ is \ $S_{22}$ \ (cf. \cite{Po}).

From now on, we shall assume that $r\ge 2$. For $r=2$, the Thom polynomial
for $I_{2,2}$ is (cf. \cite {Rim2}):
$$
S_{133}+3S_{34}\,.
$$
By virtue of Proposition \ref{Pce}, the equations from Theorem \ref{TEq}
characterizing the Thom polynomial for $I_{2,2}$ are:
\begin{equation}\label{Ia}
P(-\B_{r-1})=P(x-\fbox{$2x$}-\B_{r-1})=P(x-\fbox{$3x$}-\B_{r-1})=0\,,
\end{equation}
and
\begin{equation}\label{Ib}
P(\X_2\moins \fbox{$2x_1$}\moins \fbox{$2x_2$}\moins \B_{r-1})
=x_1x_2(x_1\moins 2x_2)(x_2\moins 2x_1)
\ R(\X_2\plus \fbox{$x_1\plus x_2$},\B_{r-1})\,.
\end{equation}
Here, without loss of generality, we assume that $x$, $x_1$, $x_2$,
and $\B_{r-1}$ are variables.
Moreover, $P(\bullet)$ denotes a symmetric function.
For the remainder of this paper, we set
\begin{equation}
\D:=\fbox{$2x_1$}+\fbox{$2x_2$}+\fbox{$x_1+x_2$}\,.
\end{equation}
Then, additionally, for variables $x_1,x_2$ and an alphabet $\B_{r-2}$,
we have the vanishing imposed by $III_{2,2}$:
\begin{equation}\label{Ic}
P(\X_2-\D -\B_{r-2})=0\,.
\end{equation}
Indeed, the singularities $\ne I_{2,2}$ with codimension
$\le \codim(I_{2,2})$ are: $A_0$, $A_1$, $A_2$, $III_{2,2}$.

For $r\ge 1$, we set
\begin{equation}
P_r(\bullet):={\cal T}^{I_{2,2}}_r(\bullet)\,.
\end{equation}

\begin{lemma}
(i) A partition appearing nontrivially in the Schur function expansion
of $P_r$ contains the rectangular partition $(r+1,r+1)$.

\smallskip

\noindent
(ii) A partition appearing nontrivially in the Schur function expansion
of $P_r$ has at most three parts.
\end{lemma}
Proof.
(i) This follows from the fact that the singularity $I_{2,2}$ belongs
to the Thom-Boardman singularity $\Sigma^2$.

\smallskip

\noindent
(ii) We can assume that $r\ge 3$. In addition to information contained in (i), we
shall use Eq. (\ref{Ic}):
$$
P_r(\X_2-\D -\B_{r-2})=0\,.
$$
By virtue of (i), we can use factorization property (\ref{Fact})
to all summands of
\begin{equation}
P_r(\X_2-\D-\B_{r-2})=\sum_I \alpha_I S_I(\X_2-\D-\B_{r-2})
\end{equation}
(we assume that $\alpha_I\ne 0$). We divide each summand of this last
polynomial by the resultant
$$
R(\X_2, \D+\B_{r-2})\,.
$$
Suppose that the resulting factor of $S_I$ is:
\begin{equation}\label{factor}
S_{p,q}(\X_2) \ S_J(-\D-\B_{r-2})\,,
\end{equation}
cf. (\ref{Fact}). Since $|I|=3r+1$, we have
\begin{equation}\label{r-1}
|J|\le r-1 \,.
\end{equation}
Now, let us assume that $I$ has more than 3 parts, that is
$J$ has more than 2 parts. This assumption (together with the
inequality (\ref{r-1})) implies that
$$
S_J(-\B_{r-2})\ne 0
$$
($\B_{r-2}$ is an alphabet of variables).
Expanding (\ref{factor}), we get among summands the following one of
largest possible degree $|J|$ in $\B_{r-2}$:
\begin{equation}\label{summ}
S_{p,q}(\X_2) \ S_J(-\B_{r-2})\ne 0\,.
\end{equation}
Take in the sum
$$
\sum_I \alpha_I S_{p,q}(\X_2) \ S_J(-\D-\B_{r-2})
$$
the (sub)sum of all the nonzero summands of the form (\ref{factor}) with
the largest possible weight of $J$.
Since Schur polynomials are independent this subsum is nonzero and moreover
it is ${\bf Z}$-linearly independent of other summands both in the sum indexed by
partitions with $\ge 3$ parts, and as well as in that indexed by partitions with 2
parts (this last sum does not depend on $\B_{r-2}$).
Hence, there is no {\bf Z}-linear combination of $S_I$'s which involve
nontrivially $I$ with more than three parts and possibly also those with 3
and 2 parts, that satisfies Eq. (\ref{Ic}). Assertion (ii) has been proved.
\qed

\smallskip

\noindent
(For example, $S_{1144}$ cannot appear in the Schur function
expansion of $P_3$ because $S_{1144}(\X_2-\D-\B_1)$ after division
by the resultant contains the summand $S_{11}(-\B_1)=S_2(\B_1)$,
which does not occur in similar expressions for $S_{55}, S_{46},
S_{244}, S_{145}$.)

\medskip

The following lemma gives a recursive description of $P_r$.
Denote by $\tau$ the linear endomorphism on the ${\bf Z}$-module of
Schur functions corresponding to partitions of length $\le 3$
that sends a Schur function $S_{i_1,i_2,i_3}$ to $S_{i_1+1,i_2+1,i_3+1}$.
Let $P^o_r$ denote the sum of those terms in the Schur function expansion
of $P_r$ which correspond to partitions of length $\le 2$.
Note that $P^o_1=S_{22}$.

\begin{lemma}\label{Lr} With this notation, for $r\ge 2$, we have the following
recursive equation:
\begin{equation}\label{lr}
P_r=P_r^o+\tau(P_{r-1})\,.
\end{equation}
\end{lemma}
\proof
Write
\begin{equation}\label{dP}
P_r=\sum_I \alpha_I S_I=\sum_J \alpha_J S_J + \sum_K \alpha_K S_K\,,
\end{equation}
where $J$ have 2 parts and $K=(k_1,k_2,k_3)$ have 3 parts
(we assume that $\alpha_I\ne 0$). We set
\begin{equation}
Q=\sum_K \alpha_K S_{k_1-1,k_2-1,k_3-1}\,,
\end{equation}
and our goal is to show that $Q=P_{r-1}$.
Since a partition $I$ appearing nontrivially in the Schur function expansion
of $P_r$ must contain the partition $(r\plus 1,r\plus 1)$, then any
partition $K$ above contains the partition $(r,r)$. Since this last partition
is not contained in the $(1,r-1)$-hook, Eqs. (\ref{Ia}) with $r$ replaced
by $r-1$ are automatically fulfilled by virtue of the vanishing property
(\ref{van}).
Note that Eq. (\ref{Ic}) is a particular case of Eq. (\ref{Ib}).
Indeed, specializing $b_{r-1}$ to \fbox{$x_1\plus x_2$} in Eq. (\ref{Ib}),
we get Eq. (\ref{Ic}). Therefore it suffices to show that
\begin{equation}
Q(\X_2-\E-\B_{r-2})=
x_1x_2(x_1\moins 2x_2)(x_2\moins 2x_1)
\ R(\X_2\plus \fbox{$x_1\plus x_2$},\B_{r-2})\,.
\end{equation}
where $\E=\fbox{$2x_1$}+\fbox{$2x_2$}$.
We apply to each summand
$$
\alpha_K S_{k_1-1,k_2-1,k_3-1}(\X_2-\E-\B_{r-2})
$$
of $Q(\X_2-\E-\B_{r-2})$
the factorization property (\ref{Fact}), and divide it by the resultant
$$
R(\X_2, \E+\B_{r-2})\,.
$$
Suppose that the resulting factor is:
\begin{equation}
\alpha_K S_{a,b}(\X_2) \ S_c(-\E-\B_{r-2})\,,
\end{equation}
where $(k_1-1,k_2-1,k_3-1)=(c,r+a,r+b)$.

Performing the same division of
$$
x_1x_2(x_1\moins 2x_2)(x_2\moins 2x_1)
\ R(\X_2\plus \fbox{$x_1\plus x_2$},\B_{r-2})
$$
we get $R(\fbox{$x_1\plus x_2$},\B_{r-2})$. Thus the wanted
equation $Q=P_{r-1}$ is equivalent to
\begin{equation}\label{show}
\sum_{a+b+c=r-2} \alpha_K S_{a,b}(\X_2) \ S_c(-\E-\B_{r-2})=
R(\fbox{$x_1\plus x_2$},\B_{r-2})\,.
\end{equation}

To prove Eq. (\ref{show}) we use Eqs. (\ref{Ib}) and (\ref{dP}) for $P_r$:
$$
\sum_I \alpha_I S_I(\X_2\moins \E \moins \B_{r-1})
=x_1x_2(x_1\moins 2x_2)(x_2\moins 2x_1)
\ R(\X_2\plus \fbox{$x_1\plus x_2$},\B_{r-1})\,.
$$
Using again the factorization property (this time w.r.t. the larger
rectangle $(r+1)^2$) and dividing both sides of the last equation by
the resultant
$$
R(\X_2, \E+\B_{r-1})\,.
$$
we get the identity
\begin{equation}\label{eqat}
\sum_{p+q+j=r-1} \alpha_I S_{p,q}(\X_2) \ S_j(-\E-\B_{r-1})=
R(\fbox{$x_1\plus x_2$},\B_{r-1}).
\end{equation}
Since
$$
S_j(-\E -\B_{r-1})=S_j(-\E-\B_{r-2})-b_{r-1}S_{j-1}(-\E -\B_{r-2})
$$
and
$$
R(\fbox{$x_1\plus x_2$},\B_{r-1})=(x_1+x_2-b_{r-1})
R(\fbox{$x_1\plus x_2$},\B_{r-2})\,,
$$
taking the coefficients of $(-b_{r-1})$ in both sides of Eq. (\ref{eqat}),
we get the wanted Eq. (\ref{show}). The lemma has been proved.
\qed

\smallskip

\noindent
(For example, writing $P_3=\alpha S_{46}+\beta S_{55}+\gamma S_{244}+\delta
S_{145}$, we get that
$$
\gamma S_1(-\E-B_1)+\delta S_1(\X_2)= R(\fbox{$x_1\plus x_2$}, \B_1)
$$
by taking the coeficients of $(-b_2)$ in both sides of
$$
\alpha S_2(\X_2)+\beta S_{11}(\X_2)+\gamma S_2(-\E-\B_2)+\delta
S_1(-\E-\B_2)S_1(X_2)
= R(\fbox{$x_1\plus x_2$}, \B_2)\,.)
$$

\smallskip

Iterating Eq. (\ref{lr}) gives
\begin{corollary}\label{Cr} With the above notation, we have
\begin{equation}\label{cr}
P_r=P^o_r+\tau(P^o_{r-1})+\tau^2(P^o_{r-2})+\cdots+\tau^{r-1}(P^o_1)\,.
\end{equation}
\end{corollary}

Of course, $P^o_r$ is uniquely determined by its value on
$\X_2$. The following result gives this value.

\begin{proposition}\label{Po}
For any $r\ge 1$, we have
\begin{equation}
P^o_r(\X_2)=(x_1x_2)^{r+1} \ S_{r-1}(\D)\,.
\end{equation}
\end{proposition}
\proof
We use induction on $r$. For $r=1,2$, the assertion holds true.
Suppose that the assertion is true for $P^o_i$ where $i<r$.
We consider the Schur function expansion of $P_r$:
\begin{equation}\label{sum}
P_r=\sum_I \alpha_I S_I\,.
\end{equation}
Fix a partition $I=(j,r+1+p,r+1+q)$ appearing nontrivially in (\ref{sum}).
Note that $j$ varies from $0$ to $r-1$ because $|I|=3r+1$.
We obtain by the factorization property (\ref{Fact}):
$$
S_I(\X_2-\D-\B_{r-2})=R \cdot S_j(-\D-\B_{r-2}) \cdot S_{p,q}(\X_2)\,.
$$
where $R=R(\X_2, \D + \B_{r-2})$. Hence, using Eq. (\ref{cr}),
we see that
\begin{equation}\label{sumj}
P_r(\X_2-\D-\B_{r-2}) = R \cdot
\Bigl(\sum_{j=0}^{r-1} S_j(-\D-\B_{r-2})
\ \frac{P_{r-j}^o(\X_2)}{(x_1x_2)^{r-j+1)}}\Bigr)\,.
\end{equation}
By the induction assumption, for positive $j\le r-1$,
$$
P_{r-j}^o(\X_2)=(x_1x_2)^{r-j+1} \ S_{r-1-j}(\D)\,.
$$
Substituting this to (\ref{sumj}), and using the vanishing (\ref{Ic}),
we obtain
\begin{equation}\label{l1}
\sum_{j=1}^{r-1} S_j(-\D-\B_{r-2}) S_{r-1-j}(\D)
+\frac{P_r^o(\X_2)}{(x_1x_2)^{r+1}}=0\,.
\end{equation}
But we also have, by a formula for addition of alphabets,
\begin{equation}\label{l2}
\sum_{j=1}^{r-1} S_j(-\D-\B_{r-2}) S_{r-1-j}(\D)+S_{r-1}(\D)
=S_{r-1}(-\B_{r-2})=0\,.
\end{equation}
Combining Eqs. (\ref{l1}) and (\ref{l2}) gives
$$
P^o_r(\X_2)=(x_1x_2)^{r+1} \ S_{r-1}(\D)\,,
$$
that is, the assertion of the induction.
The proof of the proposition is now complete.
\qed

This proposition allows us to write down the Schur function decomposition
of $P_r^0$. The coefficients of $S_{22}, S_{34}, S_{46}, S_{58},
\ldots$ (in general, $S_{i,2i-2}$), are given by the coefficients
$1, 3, 7, 15, \ldots$ in the expansion of the series:
$$
\aligned
\frac{1}{(1-z)(1-2z)}\\
=1+3z+&7z^2+15z^3+31z^4+63z^5+127z^6+\ldots\,.
\endaligned
$$
Denote these coefficients by $d_{11}, d_{21}, d_{31}
, d_{41}, \ldots$.
Moreover, we set $d_{1j}=d_{2j}=0$ for $j\ge 1$, $d_{3j}=d_{4j}=0$
for $j\ge 2$, $d_{5j}=d_{6j}=0$ for $j\ge 3$ etc.
Next, denoting by $d_{rj}$ the coefficient of $S_{r+j,2r+1-j}$ in $P_r^0$,
where $j=1,\ldots, [(r+1)/2]$, we have the recursive formula
\begin{equation}\label{rd}
d_{i+1,j}= d_{i,j-1} + d_{ij}\,.
\end{equation}
We get the following matrix:

\smallskip

$$
\begin{array}{cccccc}
d_{11} & 0 & 0 & 0 & 0 & \ldots \\
d_{21} & 0 & 0 & 0 & 0 & \ldots \\
d_{31} & d_{32} & 0 & 0 & 0 & \ldots \\
d_{41} & d_{42} & 0 & 0 & 0 & \ldots \\
d_{51} & d_{52} & d_{53} & 0 & 0 & \ldots \\
d_{61} & d_{62} & d_{63} & 0 & 0 & \ldots \\
d_{71} & d_{72} & d_{73} & d_{74} & 0 & \ldots \\
\vdots & \vdots & \vdots & \vdots & \vdots &
\end{array} \ \ \ \ \ \
= \ \ \ \ \ \
\begin{array}{cccccc}
1 & 0 & 0 & 0 & 0 & \ \ldots \\
3 & 0 & 0 & 0 & 0 & \ldots \\
7 & 3 & 0 & 0 & 0 & \ldots \\
15 & 10 & 0 & 0 & 0 & \ldots \\
31 & 25 & 10 & 0 & 0 & \ldots \\
63 & 56 & 35 & 0 & 0 & \ldots \\
127 & 119 & 91 & 35 & 0 & \ldots \\
\vdots & \vdots & \vdots & \vdots & \vdots &
\end{array}
$$

\smallskip

Summing up, we have

\begin{proposition}\label{Ppo} For $r\ge 1$, we have
\begin{equation}
P_r^o=\sum_{j=1}^{[(r+1)/2]} \ d_{rj} \ S_{r+j,2r+1-j}\,,
\end{equation}
where the $d_{ij}$'s are defined above.
\end{proposition}

We have the following values of $P_1^o, P_2^o, \ldots, P_6^o$:
$$
\aligned
S_{22}, \ 3S_{34}, \ 7S_{46}+3S_{55}, \ 15S_{58}+10&S_{67},
\ 31S_{6,10}+25S_{79}+10S_{88}\,, \\
63&S_{7,11}+56S_{8,10}+35S_{99}\,.
\endaligned
$$

Combining Proposition (\ref{Ppo}) with Eq. (\ref{cr}), we get
\begin{theorem}\label{TI22} For $r\ge 1$ the Thom polynomial
for $I_{2,2}$, with parameter $r$, equals
\begin{equation}
P_r=\sum_{k=0}^{r-1} \ \sum_{\{j\ge 1: \ k+2j\le r+1\}} \ d_{r-k,j}
\ S_{k,r+j,2r-k-j+1}\,.
\end{equation}
\end{theorem}

We have the following values of $P_1, P_2=\tau(P_1)+P_2^o, \ldots,
P_6=\tau(P_5)+P_6^o$:
$$
\aligned
&S_{22}\\
&S_{133}\plus 3S_{34}\\
&S_{244}\plus 3S_{145}\plus 7S_{46}\plus 3S_{55}\\
&S_{355}\plus 3S_{256}\plus 7S_{157}\plus 3S_{166}\plus 15S_{58}
\plus 10S_{67}\\
&S_{466}\plus 3S_{367}\plus 7S_{268}\plus 3S_{277}\plus 15S_{169}\plus
10S_{178}\plus 31S_{6,10}\plus 25S_{79}\plus 10S_{88}\\
&S_{577}\plus 3S_{489}\plus 7S_{379}\plus 3S_{388}\plus 15S_{2,7,10}
\plus 10S_{289}\plus 31S_{1,7,11}\plus 25S_{1,8,10}\plus 10S_{189}\plus \\
&63S_{7,11}\plus 56S_{8,10}\plus 35S_{99}\,.
\endaligned
$$

\section{Towards Thom polynomials for $A_i$}

The following function $F^{(i)}_r$ will be basic for computing the Thom
polynomials for $A_i$ ($i\ge 1$). We set
\begin{equation}
F^{(i)}_r(\bullet):=\sum_{J} \ S_J(\fbox{$2$}+\fbox{$3$}
+\cdots+\fbox{$i$}) S_{r-j_{i-1},\ldots,r-j_1,r+|J|}
(\bullet)\,,
\end{equation}
where the sum is over partitions $J\subset (r^{i-1})$,
and for $i=1$ we understand $F^{(1)}_r(\bullet)=S_r(\bullet)$.

\begin{example} We have
$$
F^{(2)}_r(\bullet)=\sum_{j\le r} S_j(\fbox{$2$}) S_{r-j,r+j}(\bullet)
=\sum_{j\le r} 2^j S_{r-j,r+j}(\bullet)\,,
$$
$$
F^{(3)}_r(\bullet)=\sum_{j_1\le j_2 \le r} S_{j_1,j_2}(\fbox{$2$}+\fbox{$3$})
S_{r-j_2,r-j_1,r+j_1+j_2}(\bullet)\,,
$$
$$
F^{(4)}_r(\bullet)=\sum_{j_1\le j_2\le j_3 \le r} S_{j_1,j_2,j_3}
(\fbox{$2$}+\fbox{$3$}+\fbox{$4$})
S_{r-j_3,r-j_2,r-j_1,r+j_1+j_2+j_3}(\bullet)\,,
$$
$$
F^{(i)}_1(\bullet)=\sum_{j\le i-1} \Lambda_j(\fbox{$2$}+\fbox{$3$}
+\cdots+\fbox{$i$})S_{1^{i-j-1},j+1} (\bullet)\,,
$$
where $\Lambda_j(\bullet)=(-1)^j S_j(- \ \bullet)$.
\end{example}

In the following, we shall tacitly assume that $x$, $x_1$, $x_2$, and $\B_r$
are variables (though many results remain valid without this assumption).

The following result gives the key algebraic property of $F^{(i)}_r$.
\begin{proposition}\label{FBr} \ We have
\begin{equation}\label{Br}
F^{(i)}_r(x-\B_r)= R(x+\fbox{$2x$}+\fbox{$3x$}+\cdots
+\fbox{$ix$}\,, \B_r)\,.
\end{equation}
\end{proposition}
\proof
The assertion follows from Lemma \ref{LFR} with $m=i-1$, $n=r$, and
$\A=\fbox{$2$}+\fbox{$3$}+\cdots +\fbox{$i$}$\,.
\qed

\begin{corollary}\label{CF}
Fix an integer $i\ge 1$.

\noindent
(i) For $p\le i$, we have
\begin{equation}
F^{(i)}_r(x-\B_{r-1}-\fbox{$px$})=0\,.
\end{equation}
(ii) Moreover, we have
\begin{equation}
F^{(i)}_r(x\moins \B_{r-1}\moins \fbox{$(i\plus 1)x$})=
R(x\plus \fbox{$2x$}\plus \fbox{$3x$}\plus
\cdots\plus\fbox{$ix$}\,, \B_{r-1}\plus \fbox{$(i\plus 1)x$}\, )\,.
\end{equation}
\end{corollary}
\proof Substituting in Eq. (\ref{Br}):
$$
\B_r=\B_{r-1}+\fbox{$px$}
$$
for $p\le i$, and, respectively,
$$
\B_r=\B_{r-1}+\fbox{$(i\plus 1)x$}\,,
$$
we get the assertions.\qed

\begin{theorem}(\cite{T}, \cite{Ro}) \
The polynomials \ $S_r$ \ and \ $\sum_{j\le r} 2^j S_{r-j,r+j}$
are Thom polynomials (with parameter $r$) for the singularities
$A_1$ and $A_2$.
\end{theorem}
Proof. Since only $A_0$ has smaller codimension than $A_1$,
and only $A_0$, $A_1$ are of smaller codimension
than $A_2$, the equations from Theorem \ref{TEq}
characterizing these Thom polynomials are:
\begin{equation}
P(-\B_{r-1})=0, \ \ P(x-\B_{r-1}-\fbox{$2x$})=
R(x, \B_{r-1}\plus \fbox{$2x$})
\end{equation}
for $A_1$, and
\begin{equation}
\aligned
P(-\B_{r-1})=P(x-\B_{r-1}-\fbox{$2x$})=0, \\
P(x\moins \B_{r-1}-\fbox{$3x$})&=
R(x\plus \fbox{$2x$}, \B_{r-1}\plus \fbox{$3x$}\, )
\endaligned
\end{equation}
for $A_2$. Hence the claim follows from Corollary \ref{CF}.\qed

\medskip

Since the singularities $\ne A_3$, whose codimension is $\le \codim(A_3)$
are: $A_0$, $A_1$, $A_2$ and, for $r\ge 2$, $III_{2,2}$ (cf. \cite{dPW}),
Theorem \ref{TEq} yields the following equations characterizing
${\cal T}^{A_3}$:
\begin{equation}\label{EqA3}
P(-\B_{r-1})=P(x-\B_{r-1}-\fbox{$2x$})=P(x\moins \B_{r-1}-\fbox{$3x$})=0\,,
\end{equation}
\begin{equation}
P(x- \B_{r-1}-\fbox{$4x$})=
R(x+ \fbox{$2x$}+ \fbox{$3x$}, \B_{r-1}+ \fbox{$4x$}\, )
\end{equation}
\begin{equation}
P(\X_2-\D-\B_{r-2})=0\,.
\end{equation}
By Corollary \ref{CF}, the first four equations are satisfied by
the function $F^{(3)}_r$. For $r=1$, this means that
\begin{equation}
F^{(3)}_1=S_{111}+5S_{12}+6S_3
\end{equation}
is the Thom polynomial for $A_3$. However, for $r\ge 2$, $F^{(3)}_r$
does not satisfy the last vanishing, imposed by $III_{2,2}$.
In the following we shall ``modify'' $F^{(3)}_r$ in order to obtain
the Thom polynomial for $A_3$.

\medskip

Let us discuss now $A_4$ for $r=1$ (its codimension is $4$). Then
the singularities $\ne A_4$, whose codimension is $\le \codim(A_4)$
are: $A_0$, $A_1$, $A_2$, $A_3$, $I_{2,2}$. The Thom polynomial is
\begin{equation}\label{A4,1}
{\cal T}^{A_4}=S_{1111}+9S_{112}+26S_{13}+24S_4+10S_{22}\,.
\end{equation}
This Thom polynomial was originally computed in \cite{G} via the
desingularization method. Its alternative derivation via solving
equations imposed by the above singularities was done in \cite{Rim1}).

It may be instructive for the reader to reprove here
this result using the function $F^{(4)}_1$. In this way we show
(on this relatively simple example) the method used later to more
complicated singularities. A Thom polynomial
is a sum of $j$-polynomials
(cf. Introduction) associated with the bundle morphism $df:TX \to f^*TY$.
In fact, $F^{(i)}_r$ is the $j$-polynomial part for the largest possible
$j$ (the corresponding rectangle is a row of length $r$).
Then to get the correct Thom polynomial, the function $F^{(i)}_r$
must be modified by $j$-polynomials related with smaller $j$'s.
We shall see this in the next section for the singularity $A_3$ and $r\ge 2$.
In the present case, this works as follows. We have
\begin{equation}
F^{(4)}_1=S_{1111}+9S_{112}+26S_{13}+24S_4\,.
\end{equation}
By Corollary \ref{CF}, this function satisfies the following
equations imposed by $A_0$, $A_1$, $A_2$, $A_3$, $A_4$:
\begin{equation}\label{A4}
F^{(4)}_1(0)=F^{(4)}_1(x-\fbox{$2x$})=F^{(4)}_1(x-\fbox{$3x$})=
F^{(4)}_1(x-\fbox{$4x$})=0\,,
\end{equation}

\smallskip

\begin{equation}\label{A4n}
F^{(4)}_1(x-\fbox{$5x$})= R(x+\fbox{$2x$}+\fbox{$3x$}+\fbox{$4x$},
\fbox{$5x$})\,.
\end{equation}
However, $F^{(4)}_1$ does not satisfy the vanishing imposed by $I_{2,2}$.
Namely, we have
\begin{equation}\label{F41}
F^{(4)}_1(\X_2-\fbox{$2x_1$}-\fbox{$2x_2$})
=(-10)x_1x_2(x_1-2x_2)(x_2-2x_1)\,.
\end{equation}\label{F4R}
To see this, invoke Proposition \ref{FBr}:
\begin{equation}
F^{(4)}_1(x-\B_1)=R(x+\fbox{$2x$}+\fbox{$3x$}+\fbox{$4x$},\B_1)\,.
\end{equation}
Substituting to the LHS of Eq. (\ref{F41}) $x_1=0$, we get by this
proposition
$$
F^{(4)}_1(x_2-\fbox{$2x_2$})=R(x_2+\fbox{$2x_2$}+\fbox{$3x_2$}
+\fbox{$4x_2$},\fbox{$2x_2$})=0\,,
$$
and substituting $x_1=2x_2$,
$$
\aligned
F^{(4)}_1(x_2-\fbox{$2x_1$})=R(x_2+\fbox{$2x_2$}+\fbox{$3x_2$}
+\fbox{$4x_2$},\fbox{$2x_1$})\\
=R(x_2+\fbox{$2x_2$}+\fbox{$3x_2$}+\fbox{$4x_2$},\fbox{$4x_2$})&=0\,.
\endaligned
$$
Therefore $x_1x_2(x_1-2x_2)(x_2-2x_2)$ divides this LHS. The coefficient
$-10$ results from specialization $x_1=x_2=1$.
This implies Eq. (\ref{F41}).

\smallskip

On the other hand, the Schur function $S_{22}(\bullet)$ satisfies
Eqs. (\ref{A4}), and Eq. (\ref{A4n}) with its RHS replaced by zero:
$$
S_{22}(0)=S_{22}(x-\fbox{$2x$})=\cdots = S_{22}(x-\fbox{$5x$})=0
$$
because the partition $22$ is not contained in the $(1,1)$-hook.
Moreover, we have
\begin{equation}\label{S22}
S_{22}(\X_2-\fbox{$2x_1$}-\fbox{$2x_2$})
=R(\X_2,\fbox{$2x_1$}+\fbox{$2x_2$})
=x_1x_2(x_1-2x_2)(x_2-2x_2)\,.
\end{equation}

Combining Eq. (\ref{F41}) with Eq. (\ref{S22}), the desired expression
(\ref{A4,1}) follows.

\smallskip

\begin{remark}
Porteous \cite{Po} (see also \cite{Ku}) gives a geometric account
to the function $F^{(i)}_1$. By passing with his formulas to Schur
function expansions, we should restrict ourselves to hook partitions
(more precisely: to their conjugates). We have the following recursive
formula ({\it loc. cit.}):
\begin{equation}
F^{(i)}_1=\sum_{j=1}^i \frac{(i-1)!}{(i-j)!} \ \Lambda_j \ F^{(j)}_1\,.
\end{equation}
\end{remark}

Our goal now is to give an expresion for the Thom polynomial for
$A_3$ (any $r$) as a linear combination of Schur functions.
The cases $r=1,2$ were already known in the literature (cf., e.g.,
\cite{Rim2}). In \cite{BFR}, the authors announced
a certain expression in terms of the Chern monomial basis.
Our expression is of different form (a linear combination of Schur
functions), and for the moment we do not know how to pass from it
to the one in \cite{BFR}. (A computer check for small values of $r$
shows the desired coincidence.)
The Thom polynomial for $r=1$ has been already
discussed. For $r=2$, the Thom polynomial is
\begin{equation}
S_{222}+5S_{123}+6S_{114}+19S_{24}+30S_{15}+36S_6+5S_{33}\,,
\end{equation}
and it differs from $F^{(3)}_2$ by $5S_{33}$ which is
the ``correction term'' in this case. In the following, we shall find
such a correction term for any $r$.

Define integers $e_{ij}$, for $i\ge 2$ and $j\ge 0$ in the following way.
First, $e_{20}, e_{30}, e_{40},\ldots$ are the coefficients
$5, 24, 89,\ldots$ in the development of the series:
$$
\aligned
\frac{5-6z}{(1-z)(1-2z)(1-3z)}\\
=5+24z+&89z^2+300z^3+965z^4+3024z^5+9329z^6+\ldots\,.
\endaligned
$$
Moreover, we set $e_{2j}=e_{3j}=0$ for $j\ge 1$, $e_{4j}=e_{5j}=0$
for $j\ge 2$, $e_{6j}=e_{7j}=0$ for $j\ge 3$ etc.
To define the remaining $e_{ij}$'s, we use the recursive formula
\begin{equation}\label{re}
e_{i+1,j}= e_{i,j-1} + e_{ij}\,.
\end{equation}
We now define the function $H_r(\bullet)$:
\begin{equation}\label{hr}
H_r(\bullet):=\sum_{k=0}^{r-2} \ \ \sum_{\{j\ge 0: \ k+2j \le r-2\}}
\ e_{r-k,j} \ S_{k,r+j+1,2r-k-j-1}(\bullet)\,.
\end{equation}
We state
\begin{theorem}\label{TA3}
The Thom polynomial for the singularity $A_3$, with parameter $r$, is equal
to $F^{(3)}_r+H_r$\,.
\end{theorem}
We outline a proof of the theorem in the appendix.
We shall now present some examples.
We have the following matrix $[e_{ij}]$~:

$$
\begin{array}{cccccc}
e_{20} & 0 & 0 & 0 & 0 & \ldots \\
e_{30} & 0 & 0 & 0 & 0 & \ldots \\
e_{40} & e_{41} & 0 & 0 & 0 & \ldots \\
e_{50} & e_{51} & 0 & 0 & 0 & \ldots \\
e_{60} & e_{61} & e_{62} & 0 & 0 & \ldots \\
e_{70} & e_{71} & e_{72} & 0 & 0 & \ldots \\
e_{80} & e_{81} & e_{82} & e_{83} & 0 & \ldots \\
\vdots & \vdots & \vdots & \vdots & \vdots &
\end{array} \ \ \ \ \ \
= \ \ \ \ \ \
\begin{array}{cccccc}
5 & 0 & 0 & 0 & 0 & \ \ldots \\
24 & 0 & 0 & 0 & 0 & \ldots \\
89 & 24 & 0 & 0 & 0 & \ldots \\
300 & 113 & 0 & 0 & 0 & \ldots \\
965 & 413 & 113 & 0 & 0 & \ldots \\
3024 & 1378 & 526 & 0 & 0 & \ldots \\
9329 & 4402 & 1904 & 526 & 0 & \ldots \\
\vdots & \vdots & \vdots & \vdots & \vdots &
\end{array}
$$

Consider the following matrix whose elements are two row partitions
(the symbol ``$\emptyset$" denotes the empty partition):

$$
\begin{array}{cccccc}
33 & \emptyset & \emptyset & \emptyset & \emptyset & \ldots \\
45 & \emptyset & \emptyset & \emptyset & \emptyset & \ldots \\
57 & 66 & \emptyset & \emptyset & \emptyset & \ldots \\
69 & 78 & \emptyset & \emptyset & \emptyset & \ldots \\
7,11 & 8,10 & 9,9 & \emptyset & \emptyset & \ldots \\
8,13 & 9,12 & 10,11 & \emptyset & \emptyset & \ldots \\
9,15 & 10,14 & 11,13 & 12,12 & \emptyset & \ldots \\
\vdots & \vdots & \vdots & \vdots & \vdots &
\end{array}
$$

\noindent
We use for this matrix the same ``matrix coordinates'' as for the
previous one.
Denote by $I(i,j)$ the partition occupying the $(i,j)$th place
in this matrix. So, e.g., $I(i,0)=(i+1,2i-1)$ for $i\ge 2$. For $r\ge 2$,
we set
\begin{equation}
H_r^o(\bullet):=\sum_{j\ge 0} e(r,j) \ S_{I(r,j)}(\bullet)\,.
\end{equation}
We have the following values of $H_2^0,\ldots, H_7^0$:
$$
\aligned
5S_{33}\,,
\ 24S_{45}\,,
\ 89S_{57}+24S_{66}\,,
\ 300&S_{69}+113S_{78}\,,
\ 965S_{7,11}+413S_{8,10}+113S_{99}\,, \\
3024&S_{8,13}+1378S_{9,12}+526S_{10,11}\,.
\endaligned
$$
Then we have, with the endomorphism $\tau$ defined before Lemma \ref{Lr},
\begin{equation}
H_r=H_r^o+\tau(H_{r-1})\,,
\end{equation}
or iterating,
\begin{equation}
H_r=H_r^o+\tau(H_{r-1}^o)+\tau^2(H_{r-2}^o)+\cdots+\tau^{r-2}(H_2^o)\,.
\end{equation}
We have the following values of
$H_2, H_3=\tau(H_2)+H_3^o,\ldots, H_7=\tau(H_6)+H_7^o$:
$$
\aligned
5&S_{33}\\
5&S_{144}\plus24S_{45}\\
5&S_{255}\plus24S_{156}\plus24S_{66}\plus89S_{57}\\
5&S_{366}\plus24S_{267}\plus24S_{177}\plus89S_{168}\plus113S_{78}
\plus300S_{69}\\
5&S_{477}\plus24S_{378}\plus24S_{288}\plus89S_{279}\plus113S_{189}
\plus300S_{1,7,10}\plus113S_{99}\plus413S_{8,10}\plus965S_{7,11}\\
5&S_{588}\plus24S_{489}\plus24S_{399}\plus89S_{3,8,10}\plus113S_{2,9,10}
\plus300S_{2,8,11}\plus113S_{1,10,10}\plus413S_{1,9,11}\\
\plus&965S_{1,8,12}\plus526S_{10,11}\plus1378S_{9,12}\plus3024S_{8,13}\,.
\endaligned
$$

\section{Appendix}

We present here an outline of the proof of Theorem \ref{TA3}.
The first result says that the addition of $H_r$ to $F^{(3)}_r$
is ``irrelevant''
for the conditions imposed by $A_i$, $i=0,1,2,3$.
\begin{lemma}\label{Lv} \ The function $H_r(\bullet)$ satisfies
Eqs. (\ref{EqA3}), and we have additionally
\begin{equation}
H_r(x-\B_{r-1}-\fbox{$4x$})=0\,.
\end{equation}
\end{lemma}
It follows from the lemma that it suffices to show that
\begin{equation}
(F^{(3)}_r+H_r)(\X_2-\D-\B_{r-2})=0\,,
\end{equation}
i.e. we have the vanishing imposed by $III_{2,2}$.
We look at the specialization
\begin{equation}
H_r(\X_2-\D-\B_{r-2})\,.
\end{equation}
By the factorization property (\ref{Fact}), each polynomial
\begin{equation}
S_{c,r+1+a,r+1+b}(\X_2-\D-\B_{r-2})
\end{equation}
factorizes into:
\begin{equation}
R \cdot S_c(-\D-\B_{r-2})\cdot S_{a,b}(\X_2)\,,
\end{equation}
where $R=R(\X_2,\D+\B_{r-2})$. We set
\begin{equation}
V_r(\X_2;\B_{r-2}):= \frac{H_r}{R}\,,
\end{equation}
so that Eq. (\ref{hr}) gives
\begin{equation}
V_r(\X_2;\B_{r-2})=\sum e_{r-k,j} \ S_k(-\D-\B_{r-2})
\ S_{j,r-k-j-2}(\X_2)\,,
\end{equation}
where summation is as in Eq. (\ref{hr}).

\begin{lemma}\label{Lvr} For $r\ge 2$, we have
\begin{equation}
V_r(\X_2; \B_{r-2})= \sum_{i=0}^{r-2} \ V_{r-i}(\X_2; 0) \ S_i(-\B_{r-2})\,.
\end{equation}
\end{lemma}
We look now at the specialization of $F^{(3)}_r$ imposed by $III_{2,2}$.
\begin{lemma} The polynomial $F^{(3)}_r(\X_2 - \D - \B_{r-2})$
is divisible by the resultant $R(\X_2,\D+\B_{r-2})$.
\end{lemma}
Denote by \ $-U_r(\X_2;\B_{r-2})$ \ the factor resulting from the lemma.
\begin{lemma}\label{Lur} For $r\ge 2$, we have
\begin{equation}
U_r(\X_2; \B_{r-2})= \sum_{i=0}^{r-2} \ U_{r-i}(\X_2; 0)
\ S_i(-\B_{r-2})\,.
\end{equation}
\end{lemma}
\begin{proposition}\label{Puv} \ For $r\ge 2$ we have
\begin{equation}\label{Ur0}
U_r(\X_2; 0)=
3^{r-2}\bigl(3 S_{r-2}(\X_2) - 2 S_{1,r-3}(\X_2)\bigr)=V_r(\X_2; 0)\,.
\end{equation}
\end{proposition}

Combining Lemmas \ref{Lv}, \ref{Lvr}, \ref{Lur}, and Proposition \ref{Puv},
the assertion of Theorem \ref{TA3} follows. Details will appear elsewhere.

\bigskip\medskip

\noindent
{\bf Acknowledgments} \ I thank Alain Lascoux for many helpful
discussions on the problems treated in the present paper. Though its author
is responsible for the exposition of the details, many computations
described here were done together in November 2004 and March 2005.
I am grateful to Richard Rimanyi for introducing me to his
paper \cite{Rim2}. In addition, I thank Ozer Ozturk for pointing
out several defects of the manuscript, as well as to Laszlo Feher
and Andrzej Weber for comments.

\bigskip

\noindent
{\bf Notes}

\smallskip

\noindent
1. Schur function expansions of some (other) Thom polynomials were
studied in \cite{FK} as we have been informed by Feher. Rimanyi and
Feher report that they and Komuves also observed the nonnegativity
of the Schur function expansions of Thom polynomials, cf. \cite{FK},
\cite{FR1}.

\smallskip

\noindent
2. After completion of the first version of this paper I received
the preprint \cite{FR1} containing some results on Chern monomial
expansions of Thom polynomials: an expression for the Thom series
of $I_{2,2}$, supported by a computer evidence, and an inductive
formula for Thom polynomials. Our expressions are of different form
($\bf Z$-linear combinations of Schur functions), and for the moment
we do not know how to pass from them to the ones in \cite{FR1}.

\smallskip

\noindent
3. Thom polynomials for $A_4$ and $r=3,4$ have been computed
(January 2006) via their Schur function expansion by Ozturk,
with the help of the techniques from the present paper, cf. \cite{O}.
(Note that Thom polynomial for $A_4$ and $r=2$, was computed
in \cite{Rim2}.)

\end{document}